\documentclass{article}
\usepackage{amsmath,amsfonts,amsthm,amssymb,amscd}

\renewcommand{\Im}{\mathop{\rm Im}\nolimits}

\newcommand{\p}{\partial}
\newcommand{\e}{\varepsilon}

\newcommand{\R}{{\mathbb R}}
\newcommand{\IP}{{\mathbb P}}
\newcommand{\Z}{{\mathbb Z}}
\newcommand{\E}{{\mathbb E}}
\newcommand{\T}{{\mathbb T}}

\newcommand{\eeta}{{\boldsymbol\eta}}

\newcommand{\BB}{{\cal B}}

\newcommand{\DD}{{\cal D}}
\newcommand{\EE}{{\cal E}}
\newcommand{\FF}{{\cal F}}

\newcommand{\KK}{{\cal K}}
\newcommand{\LL}{{\cal L}}

\newcommand{\NN}{{\cal N}}
\newcommand{\PP}{{\cal P}}
\newcommand{\RR}{{\cal R}}
\newcommand{\sS}{{\cal S}}

\newcommand{\XX}{{\cal X}}

\newcommand{\PPPP}{{\mathfrak P}}

\newcommand{\XXX}{{\boldsymbol{\mathit X}}}

\newcommand{\HHH}{{\boldsymbol{\mathit H}}}
\newcommand{\TTT}{{\boldsymbol{\mathit T}}}

\newcommand{\supp}{\mathop{\rm supp}\nolimits}
\newcommand{\diver}{\mathop{\rm div}\nolimits}

\theoremstyle{plain}
\newtheorem{theorem}{Theorem}[section]
\newtheorem{lemma}[theorem]{Lemma}
\newtheorem{proposition}[theorem]{Proposition}
\newtheorem{corollary}[theorem]{Corollary}
\theoremstyle{definition}
\newtheorem{definition}[theorem]{Definition}

\theoremstyle{remark}
\newtheorem{remark}[theorem]{Remark}
\newtheorem{example}[theorem]{Example}

\numberwithin{equation}{section}

\begin{document}
\author{A.~Agrachev\,$^1$, S.~Kuksin\,$^2$, A.~Sarychev\,$^3$,
A.~Shirikyan\,$^4$}
\footnotetext[1]{SISSA-ISAS, via Beirut 2--4, Trieste 34014, Italy;
Steklov Institute of Mathematics, 8~Gubkina St., 117966 Moscow,
Russia; e-mail: {\tt Agrachev@sissa.it}}
\footnotetext[2]{Department of Mathematics,
Heriot-Watt University, Edinburgh EH14 4AS, Scotland, UK;
Steklov Institute of Mathematics, 8~Gubkina St., 117966 Moscow,
Russia; e-mail: {\tt S.B.Kuksin@ma.hw.ac.uk}}
\footnotetext[3]{DiMaD, University of Florence, via C.~Lombroso 6/17,
Firenze 50134, Italy; Institute of Control Sciences (IPU), Moscow,
Russia; e-mail: {\tt Andrey.Sarychev@dmd.unifi.it}}
\footnotetext[4]{Laboratoire de Math\'ematiques, Universit\'e
de Paris--Sud~XI, B\^atiment~425, 91405 Orsay Cedex, France;
e-mail: {\tt Armen.Shirikyan@math.u-psud.fr}}
\title{On finite-dimensional projections of distributions
for solutions of randomly forced PDE's}
\date{}
\maketitle

\begin{abstract}
The paper is devoted to studying the image of probability measures on
a Hilbert space under finite-dimensional analytic maps. We establish
sufficient conditions under which the image of a measure has a density
with respect to the Lebesgue measure and continuously depends on the
map. The results obtained are applied to the 2D Navier--Stokes
equations perturbed by various random forces of low dimension.



\end{abstract}

\tableofcontents

\setcounter{section}{-1}

\section{Introduction}
\label{s0}
Let us consider the 2D Navier--Stokes equations on the torus
$\T^2\subset\R^2$: 
\begin{equation} \label{1}
\dot u+(u,\nabla)u-\nu\Delta u+\nabla p=f(t,x), \quad\diver u=0, 
\quad x\in\T^2.
\end{equation}
Here $u=(u_1,u_2)$ and $p$ are unknown velocity field and pressure,
$\nu>0$~is the viscosity, and~$f$ is an external
force. Equations~\eqref{1} are supplemented with the initial condition
\begin{equation} \label{2}
u(0)=u_0,
\end{equation}
where $u_0$ is a given function belonging to the space~$H$ of
divergence-free vector fields in~$L^2(\T^2,\R^2)$. It is well
known~\cite{Te,CF,VF} that if the right-hand side~$f$ satisfies some
mild regularity assumptions, then problem~\eqref{1}, \eqref{2} has a
unique solution~$u$ in an appropriate functional class. Our aim is to
study some qualitative properties of solutions in the situation
when~$f$ is a random process with a sufficiently non-degenerate
distribution. More precisely, let us assume that~$f$ is a stochastic
process on the positive half-line~$\R_+$ with range
in~$L^2(\T^2,\R^2)$ such that the distribution of its restriction to
any interval~$[0,T]$ is a non-degenerate decomposable measure (see
Condition~(P) in Section~\ref{s1.1}). One of the main results of this
paper says that for any finite-dimensional subspace~$F\subset H$ and
any~$t>0$ the distribution of the projection of~$u(t)$ to~$F$ has a
density with respect to the Lebesgue measure. Similar properties are
true in the case when~$f$ is a white noise in time or a sum of
independent identically distributed random forces. Furthermore, if the
random dynamical system associated with~\eqref{1} generates a Markov
process, then the above-mentioned property is valid for any stationary
distribution. These results are important since for a number of
problems in pure and applied mathematics only finitely many Fourier
components of a solution~$u$ matter; these components correspond to a
finite--dimensional subspace of~$H$.

The fact that a solution of a nonlinear SDE stirred by a degenerate
noise has a continuous density against the Lebesgue measure is very
well known. The proofs are usually based on the Malliavin calculus
(see~\cite{nualart1995} and the references therein). There have been a
few works in which various versions of the Malliavin calculus were
developed for some stochastic PDE's (for instance,
see~\cite{ocone-1988,BP-1998,DF-1998,LNP-2000,EcH,MP-2006}). In
particular, it was proved in~\cite{MP-2006} that if~$f$ is white noise
in time and sufficiently non-degenerate in the space variables, then
the distribution of the projection of solution for~\eqref{1},
\eqref{2} to any finite-dimensional subspace has a smooth density with
respect to the Lebesgue measure.

The approach of this paper is completely different and is based on the
controllability of~\eqref{1} in finite-dimensional projections and an
abstract result on the image of probability measures under analytic
mappings. We emphasise that our method does not use the Gaussian
structure of the noise, and the proof is simpler and shorter compared
to the papers quoted above.  At the same time, if the force~$f$ is
white in time, and the results of those works apply, then our
information on the density of the distribution for the projection
of~$u(t)$ is weaker than, say, that obtained in~\cite{MP-2006}.

The paper is organised as follows. In Section~\ref{s1}, we have
compiled some preliminary results on decomposable measures on Hilbert
spaces. Section~\ref{s2} contains two abstract results on the
transformation of probability measures under analytic mappings. In
Section~\ref{s3}, we apply them to the 2D Navier--Stokes equations
with different types of additive noise. Finally, in the appendix, we
prove some auxiliary results used in the main text. 

\subsection*{Notation}
Let $X$ be a Polish space, i.\,e., separable complete metric space.
We denote by~$B_X(a,R)$ the closed ball in~$X$ of radius~$R$ centred
at~$a$. If~$a$ coincides with a selected point~$\mathbf0\in X$, then
we write~$B_X(R)$. Let~$\BB(X)$ be the Borel $\sigma$-algebra on~$X$
and let~$\PP(X)$ be the family of probability measures
on~$(X,\BB(X))$. The space~$\PP(X)$ is endowed with the total
variation norm:
$$
\|\mu_1-\mu_2\|_{\rm var}
:=\sup_{\Gamma\in\BB(X)}|\mu_1(\Gamma)-\mu_2(\Gamma)|,\quad
\mu_1,\mu_2\in\PP(X).
$$
If $\mu_1,\mu_2\in\PP(X)$ and $\mu_1$ is absolutely continuous with respect
to~$\mu_2$, then we write $\mu_1\ll\mu_2$. For a random variable~$\xi$,
we denote by~$\DD(\xi)$ its distribution.

For any Banach space~$X$, we denote by $\|\cdot\|_X$ the norm in~$X$.
If $Y$ is another Banach space, then~$\LL(X,Y)$ stands for the space of
bounded linear operators from~$X$ to~$Y$. In the case $X=Y$, we shall
write~$\LL(X)$. If~$X$ is finite-dimensional, then~$\ell_X$ denotes
the Lebesgue measure on~$X$.

\smallskip
Let $J\subset\R$ be a closed interval and let $\R_+=[0,+\infty)$. We
use the following functional spaces.

\smallskip
\noindent
$C_b(X)$ is the space of bounded continuous functions $f:X\to\R$
endowed with the norm
$$
\|f\|_\infty=\sup_{x\in X}|f(x)|.
$$

\noindent
$C(J,X)$ is the space of continuous functions $u:J\to X$.

\noindent
$L^p(J,X)$ is the space of Bochner-measurable functions $u:J\to X$
such that
$$
\|u\|_{L^p(J,X)}=\biggr(\int_J\|u(t)\|_X^pdt\biggl)^{1/p}<\infty.
$$

\noindent
$L_{\rm loc}^p(\R_+,X)$ is the space of functions $u:J\to X$ whose
restriction to any finite interval $J\subset\R_+$ belongs to
$L^p(J,X)$.

\smallskip
If~$X$ is a Hilbert space and $F\subset X$ is a closed subspace, then
${\mathsf P}_F:X\to F$ denotes the orthogonal projection in~$X$
onto~$F$.

\section{Decomposable measures on Hilbert spaces}
\label{s1}
\subsection{Definitions and examples}
\label{s1.1}
Let $X$ be a separable Hilbert space with a scalar product $(\cdot,\cdot)$
and the corresponding norm $\|\cdot\|_X$. We denote by~$\BB(X)$ the Borel
$\sigma$-algebra on~$X$ and by~$\PP(X)$ the family of probability
measures on~$(X,\BB(X))$.

\begin{definition}
We shall say that a measure $\mu\in\PP(X)$ is {\it decomposable\/}
if there is an orthonormal basis $\{g_j\}\subset X$ such that
\begin{equation} \label{1.1}
\mu=\bigotimes_{j=1}^\infty\mu_j,
\end{equation}
where $\mu_j$ is the projection of~$\mu$ to the one-dimensional
space~$X_j$ generated by~$g_j$ and~$\otimes$ denotes the tensor
product of measures.
\end{definition}

\begin{example} \label{e1.2}
Let~$\mu\in\PP(X)$ be a Gaussian measure (for instance,
see~\cite{bogachev}).  It is well known that there is a vector
$a\in X$ and a self-adjoint nuclear operator $K\in\LL(X)$ such that
the characteristic function of~$\mu$ has the form
\begin{equation} \label{1.2}
\hat\mu(z)=\exp\bigl\{i(a,z)-\tfrac12(Kz,z)\bigr\}, \quad z\in X.
\end{equation}
In this case, the vector~$a$ is the mean value of~$\mu$,
$$
a=\int_Xx\mu(dx),
$$
and $K$ is the covariance operator for~$\mu$,
$$
(Kz,z)=\int_X(z,x-a)^2\mu(dx);
$$
we refer the reader to Chapter~2 in~\cite{bogachev} for more details.
Let $\{g_j\}$ be an orthonormal basis in~$X$ formed of the
eigenvectors of~$K$ and let~$\lambda_j$ be the eigenvalue of~$K$
corresponding to~$g_j$.  If a vector $z\in X$ is written in the form
$$
z=\sum_{j=1}^\infty z_jg_j,
$$
then relation~\eqref{1.2} takes the form
\begin{equation} \label{1.3}
\hat\mu(z)=\prod_{j=1}^\infty\exp\bigl\{i(a,g_j)z_j-\tfrac12\lambda_jz_j^2\bigr\}.
\end{equation}
It follows that~$\mu$ admits decomposition~\eqref{1.1} in
which~$\mu_j$ is a one-dimensional Gaussian measure with mean
value~$(a,g_j)$ and variance~$\lambda_j$.  Note also that if
$\lambda_j=0$ for some $j\ge1$, then the measure~$\mu$ is degenerate
in the sense that its support is contained in a proper affine subspace
of~$X$.
\end{example}

In what follows, we shall deal with decomposable measures possessing some
additional properties. Namely, we consider a measure $\mu\in\PP(X)$
satisfying the following condition.

\begin{itemize}
\item[\bf(P)]
The measure $\mu$ is decomposable and has a finite second moment
\begin{equation} \label{1.4}
\int_X\|x\|_X^2\mu(dx)<\infty.
\end{equation}
Moreover, every measure~$\mu_j$ in~\eqref{1.1} possesses a continuous density~$\rho_j$
with respect to the Lebesgue measure on~$X_j$.
\end{itemize}

The following lemma describes the random variables whose distribution satisfies
property~(P); its proof is obvious.

\begin{lemma} \label{l1.3}
The distribution of an $X$-valued random variable~$\xi$ satisfies
condition~{\rm(P)} if and only if~$\xi$ has the form
\begin{equation} \label{1.5}
\xi=\sum_{j=1}^\infty b_j\xi_jg_j,
\end{equation}
where $\{g_j\}$ is an orthonormal basis in~$X$, $\{\xi_j\}$ is a
sequence of scalar independent random variables such that
$\E\,\xi_j^2=1$, and $b_j>0$ are some constants satisfying the
condition
\begin{equation} \label{1.6}
\sum_{j=1}^\infty b_j^2<\infty.
\end{equation}
\end{lemma}

\begin{example} \label{e1.4}
Let $\mu\in\PP(X)$ be a Gaussian measure with a mean value $a\in X$
and a covariance operator~$K$. We claim that~$\mu$ satisfies
condition~(P) if and only if it is non-degenerate, i.e., all
eigenvalues of~$K$ are positive. Indeed, Fernique's theorem
(see~\cite{bogachev}) implies that condition~\eqref{1.4} is satisfied
for any Gaussian measure. Furthermore, it follows from~\eqref{1.3}
that the projection~$\mu_j$ of~$\mu$ to~$X_j$ is a one-dimensional
Gaussian measure with variance~$\lambda_j$. Thus, $\mu_j$ possesses a
continuous density with respect to the Lebesgue measure if and only if
$\lambda_j>0$.

For the sequel, we note that if~$\mu\in\PP(X)$ is a non-degenerate
Gaussian measure, then
\begin{equation} \label{1.7}
\mu(B)>0\quad\mbox{for any ball $B\subset X$};
\end{equation}
see Section~3.5 in~\cite{bogachev} for a proof.
\end{example}

\subsection{Support of decomposable measures}
\label{s1.2}
Let $X$ be a separable Hilbert space and let $\mu\in\PP(X)$ be a
measure possessing property~(P). Denote by~$\{g_j\}$ the orthonormal
basis in~$X$ for which representation~\eqref{1.1} holds and, for any
integer $N\ge1$, define~$X_{(N)}$ as the $N$-dimensional space
spanned by~$g_j$, $1\le j\le N$. Let
\begin{equation} \label{1.8}
\mu_{(N)}=\bigotimes_{j=1}^N\mu_j, \quad
\nu_{(N)}=\bigotimes_{j=N+1}^\infty\mu_j.
\end{equation}
If $a\in X_{(N)}$ and $R>0$, then we denote by $B_N(a,R)$ the closed ball
in~$X_{(N)}$ of radius~$R$ centred at~$a$. In the case $a=0$, we
write~$B_N(R)$.

\begin{proposition} \label{p1.5}
Let $\mu\in\PP(X)$ be a measure satisfying condition~{\rm(P)}.
Then there is a point
\begin{equation} \label{1.9}
A=\sum_{j=1}^\infty a_jg_j\in X
\end{equation}
such that\,\footnote{In what follows, we identify~$X_j$ with the real
line and regard~$\rho_j$ as a function of a real variable.}
\begin{equation} \label{1.10}
\rho_j(a_j)>0\quad\mbox{for all $j\ge1$.}
\end{equation}
Furthermore, if $A\in X$ is a point satisfying~\eqref{1.10}, then
for any integer $N\ge1$ there is $r_N>0$ such that
\begin{equation} \label{1.11}
\supp\mu\supset B_N(A_N,r_N), \quad A_N=\sum_{j=1}^Na_jg_j.
\end{equation}
\end{proposition}

\begin{proof}
To prove the existence of~$A$, denote by $A^0=\sum_j a_j^0g_j$ any
point in the support of~$\mu$.  Then $a_j^0\in\supp\mu_j$ for any
$j\ge1$. Since $\mu_j\ll\ell_{X_j}$, there is $a_j\in\R$ such that
$|a_j-a_j^0|\le j^{-1}$ and~\eqref{1.10} holds. Defining $A\in X$ by
relation~\eqref{1.9}, we obtain the required result.

\smallskip
When proving the second part of the proposition, we shall assume,
without loss of generality, that $A=0$. Let us fix any integer
$N\ge1$.  It follows from~\eqref{1.8} and~\eqref{1.10} that
\begin{equation} \label{1.12}
\supp\mu_{(N)}\supset B_N(r_N),
\end{equation}
where $r_N>0$ is sufficiently small. We claim that~\eqref{1.11} holds
for the same constant $r_N>0$. Indeed, let $x\in B_N(r_N)$ and $\e>0$.
Consider an $X$-valued random variable~$\xi$ with distribution~$\mu$.
By Lemma~\ref{l1.3}, it has the form~\eqref{1.5}. We wish to show that
\begin{equation} \label{1.13}
\IP\bigl\{\|\xi-x\|_X\le\e\bigr\}>0.
\end{equation}
To this end, define
$$
\varphi_N=\sum_{j=1}^N b_j\xi_jg_j, \quad
\psi_N=\sum_{j=N+1}^\infty b_j\xi_jg_j.
$$
It is clear that~\eqref{1.13} will be established if we prove that
\begin{equation} \label{1.14}
\IP\bigl\{\|\varphi_N-x\|_X\le\e/2\bigr\}>0, \quad
\IP\bigl\{\|\psi_N\|_X\le\e/2\bigr\}>0.
\end{equation}
The first inequality follows immediately from~\eqref{1.12}.
To prove the second, choose an integer~$M>N$ and write
\begin{equation} \label{1.15}
\|\psi_N\|_X^2=
\sum_{j=N+1}^Mb_j^2\xi_j^2+\sum_{j=M+1}^\infty b_j^2\xi_j^2=S_M+R_M.
\end{equation}
In view of~\eqref{1.6} and the relation $\E\,\xi_j^2=1$, we have
$$
\E\,R_M=\sum_{j=M+1}^\infty b_j^2\to0\quad\mbox{as $M\to\infty$}.
$$
By Chebyshev's inequality, it follows that
\begin{equation} \label{1.16}
\IP\bigl\{R_M\le\delta\bigr\}\to1\quad\mbox{as $M\to\infty$},
\end{equation}
where $\delta>0$ is an arbitrary constant. Furthermore,
inequalities~\eqref{1.10} with $a_j=0$ imply that, for any fixed $M>N$
and~$\delta>0$, we have
\begin{equation} \label{1.17}
\IP\bigl\{S_M\le\delta\bigr\}>0.
\end{equation}
Combining \eqref{1.15}~--~\eqref{1.17} and recalling that~$R_M$ and~$S_M$
are independent, we arrive at the second inequality in~\eqref{1.14}.
This completes the proof of the proposition.
\end{proof}

\subsection{A zero-one law for analytic functions}
\label{s1.3}
Let $X$ be a separable Hilbert space and let $f:X\to\R$ be a continuous
function. Recall that~$f$ is said to be {\it analytic\/} if for any $x_0\in X$
there is $\delta>0$ such that
\begin{equation} \label{1.18}
f(x)=f(x_0)+\sum_{m=1}^\infty L_m(x-x_0)\quad
\mbox{for $x\in B_X(x_0,\delta)$},
\end{equation}
where $L_m:X\to\R$ is an $m$-linear functional and the series
in~\eqref{1.18} converges regularly. The latter means that
$$
\sum_{m=1}^\infty \|L_m\|\delta^m<\infty,
$$
where $\|\cdot\|$ stands for the norm of multilinear functionals:
$$
\|L_m\|=\sup_{\|x\|_X\le1}|L_mx|.
$$
We refer the reader to~\cite{henry} or~\cite{VF} for more details.

For any continuous function $f:X\to\R$, we set
$$
\NN_f=\{x\in X: f(x)=0\}.
$$
It is well known that if $\dim X<\infty$
and~$f$ is analytic, then either $\ell_X(\NN_f)=0$ or $f\equiv0$. The
following theorem shows that a similar result is true in the
infinite-dimensional case.

\begin{theorem} \label{t1.6}
Let $f:X\to\R$ be an analytic function and let $\mu\in\PP(X)$ be a
measure possessing property~{\rm(P)}. Then
\begin{equation} \label{1.19}
\mu(\NN_f)=\mbox{$0$ or $1$}.
\end{equation}
Furthermore, if~$f$ is not identically zero, then $\mu(\NN_f)=0$.
\end{theorem}

\begin{proof}
We shall need a result from~\cite{dineen}. To formulate it, let us
introduce some definitions.

Recall that the one-dimensional spaces~$X_j$ are defined in
condition~(P) and that~$X_{(N)}$ stands for the vector space spanned
by~$X_j$, $1\le j\le N$.  Denote by~$X_{(N)}^\bot$ the orthogonal
complement of~$X_{(N)}$ in~$X$ and set
$X_{(\infty)}=\cup_NX_{(N)}$. We shall say that $A\in\BB(X)$ is a
finite zero-one $\mu$-set if
$$
\nu_{(N)}\bigl(\{y\in X_{(N)}^\bot:
\mu_{(N)}(A_N(y))=\mbox{$0$ or $1$}\}\bigr)=1
\quad\mbox{for any integer $N\ge1$},
$$
where $A_N(y)=\{x\in X_{(N)}: x+y\in A\}$, and the measures~$\mu_{(N)}$
and~$\nu_{(N)}$ are defined by~\eqref{1.8}.

According to Theorem 4 in~\cite{dineen}, if~$\mu$ is a decomposable
measure on~$X$, then for any finite zero-one $\mu$-set $A\in\BB(X)$
there is $A'\in\BB(X)$ such that
$$
A'+X_{(\infty)}=A', \quad \mu(A)=\mu(A').
$$
Applying the Kolmogorov zero-one law (see~\cite{feller}), we see
that $\mu(A')=\mbox{$0$ or $1$}$. Thus, the measure of any finite zero-one
$\mu$-set is either zero or one.

To prove~\eqref{1.19}, note that if~$f$ is analytic, then for any integer
$N\ge1$ and any $y\in X_{(N)}^\bot$, we have either
$\ell_N(\NN_f(y))=0$ or $\ell_N(X_{(N)}\setminus\NN_f(y))=0$,
where $\ell_N$ denotes the Lebesgue measure on $X_{(N)}$.
Since $\mu_{(N)}\ll\ell_N$ (see condition~(P)), we see that
$$
\mu_{(N)}(\NN_f(y))=\mbox{$0$ or $1$}\quad
\mbox{for any $y\in X_{(N)}^\bot$}.
$$
Thus, $\NN_f$ is a finite zero-one $\mu$-set, and~\eqref{1.19} follows from
what has been said above.

\smallskip
We now suppose that $f\not\equiv0$. In this case, there is $x_0\in X$
such that $f(x_0)\ne0$, and since~$\{g_j\}$ is a basis in~$X$ and~$f$
is continuous, there is no loss of generality in assuming
that $x_0\in X_{(N)}$ for some integer $N\ge1$.  In view
of~\eqref{1.19}, the required assertion will be established if we show
that
\begin{equation} \label{1.20}
\mu\bigl(\{x\in X:f(x)\ne0\}\bigr)\ne0.
\end{equation}

The restriction of~$f$ to~$X_{(N)}$ is an analytic function on a
finite-dimensional space.  It follows that
$$
f(x)\ne0\quad\mbox{for $\ell_N$-almost every $x\in X_{(N)}$},
$$
Recalling Proposition~\ref{1.5}, we can find a point $x_0\in\supp\mu$
such that $f(x_0)\ne0$. By continuity, there is $\delta>0$ such that
$$
f(x)\ne0\quad\mbox{for $x\in B(x_0,\delta)$}.
$$
This implies that $\{f(x)\ne0\}\supset B(x_0,\delta)$, and
therefore~\eqref{1.20} holds.
\end{proof}

\section{Image of measures under analytic maps}
\label{s2}

\subsection{Formulations of the results}
\label{s2.1} Let $(H,d)$ be a metric space and let~$X$ and~$F$ be
finite-dimensional vector spaces. Consider a continuous operator
$f:H\times X\to F$. For any probability measure $\mu\in\PP(X)$ and
any $u\in H$, denote by~$f_*(u,\mu)$ the image of~$\mu$
under~$f(u,\cdot)$.

\begin{theorem} \label{t2.1}
Suppose that, for any $u\in H$, the function~$f(u,\cdot)$ is analytic and
the interior of the set~$f(u,X)$ is non-empty. Let $\mu\in\PP(X)$ be a
measure possessing a continuous density~$\rho(x)$ with respect to the
Lebesgue measure on~$X$. Then the following assertions hold.
\begin{itemize}
\item[\bf(i)]
For any $u_0\in H$, the measure~$f_*(u_0,\mu)$ is absolutely continuous
with respect to~$\ell_F$.
\item[\bf(ii)]
The function~$f_*(\cdot,\mu)$ from~$H$ to the space~$\PP(F)$ endowed
with the total variation norm is continuous.
\end{itemize}
\end{theorem}

Our next goal is to study the case in which~$X$ is an infinite-dimensional space.
More precisely, suppose that~$X$ is a separable Hilbert space and~$F$ is a
finite-dimensional vector space. Recall that condition (P) is introduced
in Section~\ref{s1.1}.

\begin{theorem} \label{t2.2}
Let $f:H\times X\to F$ be a continuous function such that~$f(u,\cdot)$
is analytic for any $u\in H$ and the derivative $D_xf(u,x)$ is
continuous with respect to $(u,x)$. Suppose that, for any $u\in H$,
there is a ball~$B_u$ in a finite-dimensional subspace~$X_u\subset X$
such that the interior of the set $f(u,B_u)$ is non-empty. Then for
any measure $\mu\in\PP(X)$ satisfying condition~{\rm(P)}
statements~{\rm(i)} and~{\rm(ii)} of Theorem~\ref{t2.1} take place.
\end{theorem}

\subsection{Proof of Theorem~\ref{t2.1}}
\label{s2.2}
Let us fix any point $u_0\in H$ and denote by~$D_xf(u_0,x)$ the
derivative (Jacobian) of the map $f(u_0,\cdot):X\to F$ at the
point~$x\in X$.  We first show that the matrix~$D_xf(u_0,x)$ has a
minor~$m(x)$ of the size~$\dim F$ such that
\begin{equation} \label{2.1}
m(x)\ne0\quad\mbox{for $\mu$-almost every $x\in X$}.
\end{equation}
To this end, recall that a point $x\in X$ is said to be {\it
regular\/} for~$f(u_0,\cdot)$ if the rank of~$D_xf(u_0,x)$ is maximal.
Any point that is not regular is said to be {\it singular\/}. In view
of Sard's theorem (see~\cite{Sternberg}), the image under the smooth
function~$f(u_0,\cdot)$ of the set of its singular points has zero
Lebesgue measure. Since the interior of~$f(u_0,X)$ is non-empty, we
conclude that~$f(u_0,\cdot)$ has a regular point $x_0\in X$.

Let $m(x)$ be the minor of~$D_xf(u_0,x)$ that is non-zero at~$x_0$.
Since~$m(x)$ is analytic, we see that $m(x)\ne0$ almost everywhere
with respect to the Lebesgue measure~$\ell_X$. Since~$\mu$ is
absolutely continuous with respect to~$\ell_X$, we conclude
that~\eqref{2.1} holds.  

For any $\e>0$, we denote
$$
X^\e=\{x\in X:|x|\le\e^{-1}, |m(x)|\ge\e\}.
$$
Then
$$
\nu_\e:=\mu(X\setminus X^\e)\to0\quad\text{as}\quad\e\to0.
$$
Let us take any $x\in X^\e$ and write it as $x=(x_1,x_2)$, where~$x_1$
denotes the variables entering the minor~$m(x)$.  Accordingly, the
space~$X$ can be represented as a direct product $X=X_1\times
X_2$. Applying the implicit function theorem, we can find open balls
$V_1\subset X_1$ and $V_2\subset X_2$ such that for any $x_2\in V_2$
and $u\in B^\e=B_X(u_0,r_\e),$ the map $f(u,\cdot,x_2)$ is a
diffeomorphism of the domain~$V_1$ onto its image~$W(u,x_2)$. Here
$r_\e>0$ is a constant that goes to zero with~$\e$.  Accordingly, we
can write~$x_1$ in terms of $u\in B^\e$, $x_2\in V_2$ and
$y=f(u,x_1,x_2)\in W(u,x_2)$:
$$
x_1=g(u,y,x_2).
$$
The sets $V=V_1\times V_2$ corresponding to various $x\in X^\e$ form
an open cover of the compact set~$X^\e$. Let us find a finite
sub-cover~$\{V^j\}$. We denote by~$\{\varphi_j(x)\}$ a continuous
partition of unity on~$X^\e$ subordinate to~$\{V^j\}$. That is, $\varphi_j\ge0$,
$\supp\varphi_j\subset V^j$, and $(\sum\varphi_j)(x)=1$ for $x\in X^\e$. Let
$\mu_j=\varphi_j\mu$.  Then
\begin{equation} \label{n0}
\mu_j=\varphi_j(x)\rho(x)\,dx_1\,dx_2= 
\varphi_j(x_1,x_2)\rho(x_1,x_2)|m(u,x_1,x_2)|^{-1}\,dy\,dx_2\,,
\end{equation}
where $m(u,x)$ is the minor of $D_xf(u,x)$ corresponding to~$x_1$, and
$x_1=g(u,y,x_2)$ on the right-hand side of~\eqref{n0}.  
Hence, $f_*(u,\mu_j)=g_j(u,y)\,dy$, where
$$
g_j(u,y)=\int\tilde\varphi_j(y,x_2)\tilde\rho(y,x_2) 
|\widetilde m(u,y,x_2)|^{-1}\,dx_2\,
$$
with $\tilde\varphi_j(y,x_2)=\varphi_j(g(u,y,x_2),x_2)$, etc. 
Let us denote $\mu_\e=(\sum\varphi_j)\mu$. Then
$\|\mu-\mu_\e\|_{\rm var}\le\nu_\e$. We have
\begin{equation} \label{n1}
f_*(u,\mu_\e)=g_\e(u,y)\,dy\,,
\end{equation}
where $g_\e(u,y)=\sum g_j(u,y)$ is a continuous function of 
$u\in B^\e$ and~$y$. Clearly,
\begin{equation}\label{n2}
\|f_*(u,\mu_\e)-f_*(u,\mu)\|_{\rm var}\le\nu_\e\quad
\mbox{for any $u\in B^\e$}.
\end{equation}

Relations~\eqref{n1} and~\eqref{n2} with $u=u_0$ and $\e\to0$ imply
assertion~(i).  Indeed, for any measurable set $Q\subset F$ with zero
Lebesgue measure, we have
$$
f_*(u,\mu)(Q)=(f_*(u,\mu)-f_*(u,\mu_\e))(Q)+f_*(u,\mu^\e)(Q),
$$
so $f_*(u,\mu)(Q)\le \nu_\e$ for each~$\e$. Hence, $f_*(u,\mu)(Q)=0$.

To prove~(ii), we fix any~$\gamma>0$ and choose~$\e>0$ such that
$\nu_\e<\frac13\,\gamma$.  Due to~\eqref{n2}, for~$u\in B^\e$ we have
\begin{equation*} \begin{split}
\|f_*(u,\mu)-f_*(u_0,\mu)\|_{\rm var} & \le 2\nu_\e+
\|f_*(u,\mu_\e)-f_*(u_0,\mu_\e)\|_{\rm var}\\
&\le2\nu_\e+ \frac12\,
\int|g_\e(u,y)-g_\e(u_0,y)|\,dy\le\gamma,
\end{split}\end{equation*}
if~$u$ is sufficiently close to~$u_0$. Since~$u_0$ is an arbitrary
point, what has been said implies~(ii).

\begin{remark} \label{r2.3}
Analysing the proof given above, one easily sees that, instead of
assuming the existence of interior points for the set $f(u,X)$, we
could require that the function $f(u,\cdot)$ should have at least one
regular point for any $u\in H$ (cf. the beginning of the proof).
\end{remark}

\subsection{Proof of Theorem~\ref{t2.2}}
\label{s2.3}
{\it Step 1}.
Let us fix any $u_0\in H$ and show that there is finite-dimensional
subspace $X_1\subset X$ spanned by some vectors of the basis~$\{g_j\}$
such that $\dim X_1=\dim F$ and~\eqref{2.1} holds, where~$m(x)$ denotes
the determinant of the matrix for the restriction of~$D_xf(u_0,x)$
to~$X_1$.  Indeed, by the hypothesis, there is a ball~$B_{u_0}$ in a
finite-dimensional subspace $X_{u_0}\subset X$ such that the interior
of~$f(u_0,B_{u_0})$ is non-empty. Since~$\{g_j\}$ is a basis in~$X$, for
any~$\delta>0$ there is a finite-dimensional subspace~$Y\subset X$
spanned by some vectors of~$\{g_j\}$ such that
\begin{equation} \label{2.12}
\|{\mathsf P}_Yy-y\|_X\le\delta\quad\mbox{for any $y\in B_{u_0}$},
\end{equation}
where ${\mathsf P}_Y:X\to X$ is the orthogonal projection in~$X$ onto
the subspace~$Y$. Now note that~$f(u_0,\cdot)$ is continuous
and~$B_{u_0}$ is compact. Therefore for any~$\e>0$ we can
find~$\delta>0$ such that
\begin{equation} \label{2.13}
\|f(u_0,z)-f(u_0,y)\|_F\le\e\quad
\mbox{for $y\in B_{u_0}$, $\|z-y\|_X\le\delta$}.
\end{equation}
Combining~\eqref{2.12} and~\eqref{2.13}, we see that for any~$\e>0$
there is a finite-dimensional subspace~$Y\subset X$ spanned by some
vectors of~$\{g_j\}$ such that
$$
\|f(u_0,{\mathsf P}_Yy)-f(u_0,y)\|_F\le\e\quad
\mbox{for $y\in B_{u_0}$}.
$$
Choosing~$\e>0$ sufficiently small and applying Proposition~\ref{p4.1}
of the Appendix (see Section~\ref{s4}), we conclude that the interior
of the set~$f(u_0,Y)$ is non-empty. Since $\dim Y<\infty$, Sard's
theorem implies that the function $f(u_0,\cdot):Y\to F$ has at least
one regular point~$x_0\in Y$. Let us denote by~$X_1$ a subspace
spanned by some vectors of the basis~$\{g_j\}$ such that the
restriction of~$D_xf(u_0,x_0)$ to~$X_1$ is an isomorphism from~$X_1$
onto~$F$. Then $m(x_0)\ne0$. Since~$m(x)$ is an analytic function, the
required result follows from Theorem~\ref{t1.6}.

\smallskip
{\it Step 2}.
We now repeat the argument used in the proof of
Theorem~\ref{t2.1}. Let us represent~$X$ as the direct product
$X=X_1\times X_2$, where~$X_1$ is constructed in Step~1 and~$X_2$
denotes the orthogonal complement of~$X_1$ in~$X$. Denote
by~$\lambda_1$ and~$\lambda_2$ the projections of~$\mu$ to the
subspaces~$X_1$ and~$X_2$, respectively, and by~$\rho(x_1)$ the
density of~$\lambda_1$ with respect to~$\ell_{X_1}$. For any~$\e>0$,
let us choose a compact set $X^\e\subset X$ such that
$$
\mbox{$\|x\|_X\le\e^{-1}$, $|m(x)|\ge\e$ for $x\in X^\e$},
\qquad \nu_\e:=\mu(X\setminus X^\e)\le\e.
$$
As in the proof of Theorem~\ref{t2.1}, we can find a constant~$r_\e>0$
going to zero with~$\e$ and a finite cover~$\{V^j=V_1^j\times V_2^j\}$
of the compact set~$X^\e$ such that~$V_1^j\subset X_1$
and~$V_2^j\subset X_2$ are balls, and for any $x_2\in V_2$
and $u\in B^\e=B_X(u_0,r_\e),$ the map $f(u,\cdot,x_2)$ is a
diffeomorphism of the domain~$V_1$ onto its image~$W(u,x_2)$. We
denote by $x_1=g(u,\cdot,x_2)$ the inverse function
of~$f(u,\cdot,x_2)$. Let~$\{\varphi_j(x)\}$ be a continuous
partition of unity on~$X^\e$ subordinate to~$\{V^j\}$ and let
$\mu_j=\varphi_j\mu$.  Then we have (cf.~\eqref{n0})
$$
\mu_j=\varphi_j(x)\rho(x_1)\,dx_1\,\lambda(dx_2)
=\varphi_j(x_1,x_2)\rho(x_1)|m(u,x_1,x_2)|^{-1}\,dy\,\lambda(dx_2)\,,
$$
where $m(u,x)$ denotes the determinant of the restriction of
$D_xf(u,x)$ to~$X_1$, and $x_1=g(u,y,x_2)$ on the right-hand side of
the formula. The rest of the proof is literally the same as that of
Theorem~\ref{t2.1}, and therefore we omit it.

\begin{remark} \label{r2.4}
The proof given above implies that the claim of Theorem~\ref{t2.2}
remains true if we replace the condition of existence of interior
points for the set~$f(u,B_u)$ by the following one: for any $u\in H$,
there is a point $x_u\in X$ and a finite-dimensional subspace
$X_u\subset X$ such that the restriction of $D_xf(u,x_u)$ to~$X_u$ is an
isomorphism from~$X_u$ onto~$F$ (cf.~Step~1 of the proof).
\end{remark}

\section{Applications}
\label{s3}
Throughout this section, we use the standard functional spaces~$H$
and~$V$ arising in the theory of Navier--Stokes equations; they are
defined in Subsection~\ref{s3.1}. We shall also use the spaces
$$
\XX=C(\R_+,H)\cap L_{\rm loc}^2(\R_+,V), \quad
\XX_T=C(J_T,H)\cap L^2(J_T,V),
$$
where $T>0$ and $J_T=[0,T]$.

\subsection{Navier--Stokes equations perturbed by a
time-discrete random force}
\label{s3.1}
Let us consider the 2D Navier--Stokes (NS) system on the
torus~$\T^2=\R^2/2\pi\Z^2$. Define the spaces
$$
H=\{u\in L^2(\T^2,\R^2): \diver u=0\mbox{ on $\T^2$}\}, \quad
V=H^1(\T^2,\R^2)\cap H,
$$
endowed with natural norms. Here $H^s(\T^2,\R^2)$ denotes the space of
vector functions $(u_1,u_2)$ whose components belong to the Sobolev
space of order~$s$. Let $\Pi:L^2(\T^2,\R^2)\to H$ be the orthogonal
projection in~$L^2(\T^2,\R^2)$ onto~$H$. After applying the
projection~$\Pi$, the NS system reduces to the following evolution
equation in~$H$:
\begin{equation} \label{3.1}
\dot u+\nu Lu+B(u,u)=g(t,x),
\end{equation}
where $\nu>0$ is the viscosity, $L=-\Pi\Delta$, and
$B(u,v)=\Pi((u,\nabla)v)$.  In this subsection, we assume that
\begin{equation} \label{3.2}
g(t,x)=\sum_{k=1}^\infty I_{k,T}(t)\eta_k(x),
\end{equation}
where $T>0$~is a parameter, $I_{k,T}(t)$~is the indicator function of
the time interval $[(k-1)T,kT)$, and~$\{\eta_k\}$ is a sequence of
$H$-valued i.i.d.~random variables defined on a probability space
$(\Omega,\FF,\IP)$. Standard theorems on well-posedness of the 2D NS
system (e.g., see~\cite{CF}) imply that, for almost every
$\omega\in\Omega$, problem~\eqref{3.1}, \eqref{3.2} has a unique
solution $u\in\XX$ that satisfies the initial condition
\begin{equation} \label{3.3}
u(0)=u_0,
\end{equation}
where $u_0\in H$ is an arbitrary function. We shall denote by
$S_t:H\to H$ the random operator that takes~$u_0$ to~$u(t)$. Our aim
is to study the distribution for projections of the random
variables~$S_{kT}(u_0)$ to finite-dimensional subspaces of~$H$.

Let~$\{e_j\}$ be a complete set of eigenfunctions for~$L$ indexed in
an increasing order of the corresponding eigenvalues~$\alpha_j$ and
let~$H_N$ be the vector space spanned by~$e_j$, $j=1,\dots,N$. We
shall assume that the i.i.d.~random variables~$\eta_k$ satisfy the
following condition.
\begin{itemize}
\item[\bf(D)]
The random variables~$\eta_k$ have the form
\begin{equation} \label{3.4}
\eta_k=\sum_{j=1}^\infty b_j\xi_{jk}e_j,
\end{equation}
where $b_j\ge0$ are some constants such that
\begin{equation} \label{3.5}
\sum_{j=1}^\infty b_j^2<\infty,
\end{equation}
and $\xi_{jk}$ are independent scalar random variables whose
distribution~$\pi_j$ possesses a density with respect to the Lebesgue
measure, and $\supp\pi_j\ni0$ for any $j\ge0$.
\end{itemize}

\begin{theorem} \label{t3.1}
Suppose that condition~{\rm(D)} is fulfilled. Then there is an
integer~$N\ge1$ not depending on~$\nu$ and~$\{\eta_k\}$ such that the
following two statements hold, provided that
\begin{equation} \label{3.6}
b_j\ne0\quad\mbox{for $j=1,\dots,N$}.
\end{equation}
\begin{itemize}
\item[\bf(i)]
For any constant $\nu>0$ and any finite-dimensional subspace
$F\subset H$, there is a discrete subset
$\TTT=\TTT(\nu,F)\subset\R_+\setminus\{0\}$ such that if $T\notin\TTT$
and~$R>0$, then for any $u\in B_H(R)$ and an appropriate integer
$k=k(\nu,F,R)\ge1$ the distribution of ${\mathsf P}_FS_{kT}(u)$
possesses a density with respect to~$\ell_F$.
\item[\bf(ii)]
Let us set $\lambda_{u}(t)=\DD({\mathsf P}_FS_t(u))$. Then $\lambda_{u}(kT)$
continuously depends on $u\in B_H(R)$ in the total variation norm.
\end{itemize}
\end{theorem}

\begin{remark} \label{r3.2}
It is possible to give a more precise description of the integer~$N$
in~\eqref{3.6}. Namely, it is the minimal integer~$N\ge1$ such that
the vectors $e_1,\dots,e_N$ form a saturating set (see Section~\ref{s4.3}
in the Appendix for a definition of a saturating set). In particular,
if the eigenfunctions of~$L$ are indexed in a suitable way, then one can
take~$N=6$.
\end{remark}

\begin{proof}[Proof of Theorem~\ref{t3.1}]
{\it Step~1}.
We wish to apply Theorem~\ref{t2.2} and Remark~\ref{r2.4}. Let~$H_0\subset H$
be the subspace spanned by those vectors~$e_j$ for which
$b_j\ne0$ and let~$\XXX_k$ be the direct product of~$k$ copies of~$H_0$.
We fix a finite-dimensional subspace~$F\subset H$ and consider the
operator
$$
f_k:\R_+^*\times H\times\XXX_k\to F,\quad
(T,u_0,\eta_1,\dots,\eta_k)\mapsto{\mathsf P}_F u(kT),
$$
where $\R_+^*=\R_+\setminus\{0\}$ and
$u(t)$ denotes the solution of problem~\eqref{3.1}~--~\eqref{3.3} in
which~$\{\eta_k\}$ is regarded as a sequence of deterministic functions in~$H_0$.
In view of Proposition~\ref{p4.2}, the operator~$f_k$ is analytic
on $\R_+^*\times H\times\XXX_k$. Furthermore, if condition~(D) is fulfilled,
then for any integer~$k\ge1$ the distribution of the $\XXX_k$-valued random variable
$\eeta_k=(\eta_1,\dots,\eta_k)$ satisfies property~(P).
Suppose we have established the existence of a discrete subset $\TTT\subset\R_+^*$
possessing the following property:
\begin{itemize}
\item[\bf(C)]
for any~$R>0$ there is an integer~$k\ge1$ such that if $T\notin\TTT$ and
$u_0\in B_H(R)$, then the derivative
\begin{equation} \label{3.7}
(D_{\eeta_k}f_k)(T,u_0,\eeta_k):\XXX_k\to F
\end{equation}
is surjective for at least one point $\eeta_k=(\eta_1,\dots,\eta_k)\in\XXX_k$.
\end{itemize}
In this case, statements~(i) and~(ii) of the theorem are straightforward
consequences of Theorem~\ref{t2.1} and Remark~\ref{r2.3} in which~$X$ and~$H$ are
replaced by~$\XXX_k$ and~$B_H(R)$.

\smallskip
{\it Step~2}.  To prove~(C), we first assume that $u_0=0$. Let us
denote by $\RR_1:L^2(J_1,H)\to H$ the operator that takes each
function $g\in L^2(J_1,H)$ to~$u(1,x)$, where $u\in\XX_1$ is the
solution of~\eqref{3.1}, \eqref{3.3} with~$u_0=0$. By
Proposition~\ref{p4.5}, there is an integer~$N\ge1$ such that the
Navier--Stokes system~\eqref{3.1} with $g\in L^2(J_1,H_N)$ is solidly
controllable in time~$1$ for the projection to~$F$. (See
Definition~\ref{d4.4} for the concept of solid controllability.) In
particular, there is a compact subset $K\subset L^2(J_1,H)$ and a
constant~$\e>0$ such that $\varPhi(K)\supset B_F(1)$ for any
continuous map satisfying the inequality
\begin{equation} \label{3.8}
\sup_{g\in K}\|\varPhi(g)-{\mathsf P}_F\RR_1(g)\|_F\le\e.
\end{equation}
For any integer~$m\ge1$, denote by~$Y_m\subset L^2(J_1,H_N)$ the
subspace of functions that are constant on any interval of the form
$\bigl[\frac{l-1}{m},\frac{l}{m}\bigr)$, $1\le l\le m$.
It is clear that $\cup_mY_m$ is dense in~$L^2(J_1,H_N)$. Since~$K$ is compact,
for any~$\delta>0$ we can find an integer~$m\ge1$ such that
$$
\sup_{g\in K}\|g-{\mathsf P}_{Y_m}g\|_V\le\delta.
$$
It follows from the continuity of~$\RR_1$ that for any~$\e>0$ there is an
integer~$m\ge1$ such that~\eqref{3.8} is satisfied for
$\varPhi(g)={\mathsf P}_F\RR_1({\mathsf P}_{Y_m}g)$. This implies that
\begin{equation} \label{3.9}
{\mathsf P}_F\RR_1(Y_m)\supset B_F(1).
\end{equation}
Now note that if we denote by~$I_{l,m}(t)$ the indicator function of
the interval $\bigl[\frac{l-1}{m},\frac{l}{m}\bigr)$ and identify the
function
$$
g(t,x)=\sum_{k=1}^{m} I_{l,m}(t)\eta_k(x)\in Y_m
$$
with the vector $\eeta=(\eta_1,\dots,\eta_m)\in\XXX_m$, then we can write
\begin{equation} \label{3.10}
f_m(m^{-1},0,\eeta_m)={\mathsf P}_F\RR_1(g).
\end{equation}
Combining~\eqref{3.9} and~\eqref{3.10}, we conclude that there is a
finite-dimensional subspace~$Y\subset Y_m$ such that
$$
f_m(m^{-1},0,Y)\supset B_F(1).
$$
Sard's theorem now implies that the derivative
$$
(D_{\eeta_m}f_m)(m^{-1},0,\eeta_m):\XXX_m\to F
$$
is surjective for at least one point $\eeta_m^0\in\XXX_m$. Since
$(D_{\eeta_m}f_m)(T,0,\eeta_m^0)$ is an analytic function with respect
to~$T>0$, we see that there is a discrete set $\TTT\subset\R_+^*$ such
that
\begin{equation} \label{3.11}
(D_{\eeta_m}f_m)(T,0,\eeta_m^0):\XXX_m\to F
\quad\mbox{is surjective for any $T\notin\TTT$}.
\end{equation}

\smallskip
{\it Step~3}.
We can now verify property~(C). Let us fix any $T\notin\TTT$ and~$R>0$.
We claim that if an integer $l\ge1$ is sufficiently large and $k=l+m$, then the
linear operator~\eqref{3.7} is surjective for $\eeta_k=(0,\dots,0,\eeta_m^0)$
and any $u_0\in B_H(R)$. Indeed, the definition of~$f_k$ implies that
$$
f_k(T,u_0,\eeta_k^0)=f_m(T,u(lT),\eeta_m^0),
$$
where~$u(t)$ is the solution of~\eqref{3.1}, \eqref{3.3} with $g\equiv0$.
It follows that
\begin{equation} \label{3.12}
\Im\{(D_{\eeta_k}f_k)(T,u_0,\eeta_k^0)\}\supset
\Im\{(D_{\eeta_m}f_m)(T,u(lT),\eeta_m^0)\}
\end{equation}
where $\Im\{A\}$ denotes the image of a linear operator~$A$. In view of~\eqref{3.11}
and the continuity of~$D_{\eeta_m}f_m$, we can find~$r>0$ such that
\begin{equation} \label{3.13}
\Im\{(D_{\eeta_m}f_m)(T,v,\eeta_m^0)\}=F
\quad\mbox{for any $v\in B_H(r)$},
\end{equation}
Furthermore, the dissipation property of the homogeneous Navier--Stokes system
implies that
\begin{equation} \label{3.14}
u(lT)\in B_H(r)\quad\mbox{for any $u_0\in B_H(R)$}.
\end{equation}
where $l=l(R,T)\ge1$ is sufficiently
large. Combining~\eqref{3.12}~--~\eqref{3.14}, we arrive at the
required result.
\end{proof}

\begin{remark}
Analysing the proof given above, it is possible to establish the
following property, which shows that the ``bad'' subset~$\TTT$
constructed in Theorem~\ref{t3.1} cannot accumulate to zero if we
allow the integer~$k$ to depend on~$T$:
\begin{itemize}
\item
For any constant $\nu>0$ and any finite-dimensional subspace $F\subset
H$ there is~$T_0>0$ such that if $T\in(0,T_0]$ and~$R>0$, then for any
$u\in B_H(R)$ and an appropriate integer $k=k(\nu,F,R,T)\ge1$ the
distribution of ${\mathsf P}_FS_{kT}(u)$ possesses a density with
respect to~$\ell_F$ and continuously depends on~$u$ in the total
variation norm.
\end{itemize}
\end{remark}

We now study stationary solutions of~\eqref{3.1},
\eqref{3.2}. Since~$\{\eta_k\}$ are i.i.d.~random variables in~$H$,
for any deterministic initial function~$u_0$ the sequence
$\{u(kT),k\ge0\}$ is a Markov chain. Thus, the set of all
solutions restricted to the times~$kT$ form a Markov family
in~$H$. Let $P_k(u,\Gamma)$ be the corresponding transition function
and let~$\PPPP_k^*$ be the Markov operator associated
with~$P_k$. Using standard a priori estimates for solutions of the
Navier--Stokes system and applying the Bogolyubov--Krylov argument
(e.g., see~\cite{Has1980}), one can show that~$\PPPP_k^*$ has at least
one stationary distribution~$\mu$:
\begin{equation} \label{stat}
\PPPP_k^*\mu=\mu\quad\mbox{for all $k\ge1$}.
\end{equation} 
A simple consequence of Theorem~\ref{t3.1} is the following result:

\begin{corollary} \label{c3.4}
Suppose that the conditions of Theorem~\ref{t3.1} are fulfilled. Let
$F\subset H$ be a finite-dimensional subspace and let
$T\notin\TTT(\nu,F)$. Then for any stationary measure~$\mu$ the
projection~${\mathsf P}_F\mu$ possesses a density with respect to the
Lebesgue measure on~$F$.
\end{corollary}

\begin{proof}
Let $\Gamma\subset F$ be a Borel set of zero Lebesgue measure and 
let~$\e>0$. Choose a constant~$R>0$ so large that 
\begin{equation} \label{i1}
\mu(B_H(R))\ge1-\e.
\end{equation} 
By assertion~(i) of Theorem~\ref{t3.1}, we can find an integer~$k\ge1$
such that
\begin{equation} \label{i2}
({\mathsf P}_FP_k)(u,\Gamma)=0\quad\mbox{for all $u\in B_H(R)$}.
\end{equation} 
Combining~\eqref{stat}~--~\eqref{i2}, we derive
$$
{\mathsf P}_F\mu(\Gamma)=\int_H({\mathsf P}_FP_k)(u,\Gamma)\mu(du)
=\int_{B_H^c(R)}({\mathsf P}_FP_k)(u,\Gamma)\mu(du)\le\mu(B_H^c(R))\le\e,
$$
where $B_H^c(R)$ denotes the complement of~$B_H(R)$.  Since~$\e>0$ was
arbitrary, we conclude that ${\mathsf P}_F\mu(\Gamma)=0$. This
completes the proof of the corollary.
\end{proof}

\subsection{Navier--Stokes equations with a non-degenerate
random perturbation}
\label{s3.2}
In this subsection, we consider the NS system~\eqref{3.1} with a right-hand side
of the form
\begin{equation} \label{3.15}
g(t,x)=h(t,x)+\eta(t,x).
\end{equation}
We assume that~$h\in L_{\rm loc}^2(\R_+,H)$ is a deterministic function
and~$\eta$ is a random process whose trajectories belongs
to~$L_{\rm loc}^2(\R_+,H_0)$, where $H_0\subset H$ is
a closed subspace. Denote by~$\mu_T$, $T>0$, the distribution of the restriction
of~$\eta$ to the interval~$J_T$ and by $S_t:H\to H$ the random operator that
takes each~$u_0\in H$ to~$u(t)$, where $u\in\XX$ is the solution of~\eqref{3.1},
\eqref{3.3}, \eqref{3.15}.

\begin{theorem} \label{t3.2}
There is an integer~$N\ge1$ not depending on~$\nu>0$ such that if
$h\in L_{\rm loc}^2(\R_+,H)$ is a given function, $T>0$ is a constant,
$H_0\subset H$ is a subspace containing~$H_N$, and~$\mu_T$ is a
decomposable measure on~$L^2(J_T,H_0)$ satisfying property~{\rm(P)},
then for any $u_0\in H$ and any positive $t\in J_T$ the following
statements take place.\,\footnote{See Remark~\ref{r3.2} for a more precise description of~$N$.}
\begin{itemize}
\item[\bf(i)]
Let $F\subset H$ be a finite-dimensional subspace and let $\lambda_u(t)$
be the distribution of~${\mathsf P}_FS_t(u)$. Then
$\lambda_u(t)\ll\ell_F$.
\item[\bf(ii)]
The measure $\lambda_u(t)$ continuously depends on $u\in H$ in the total
variation norm.
\end{itemize}
\end{theorem}

\begin{proof}
Both assertions are straightforward consequences of Theorem~\ref{t2.2}
and Propositions~\ref{p4.2} and~\ref{p4.5}. Indeed, let $f_t:H\times
L^2(J_T,H)\to F$ be the operator that takes each pair~$(u_0,\eta)$
to~${\mathsf P}_Fu(t)$, where~$u\in\XX_T$ is the solution of
problem~\eqref{3.1}, \eqref{3.3}, \eqref{3.15} with deterministic
functions~$h$ and~$\eta$. By assumption, the distribution~$\mu_T$ of
the restriction of~$\eta$ to~$J_T$ satisfies condition~(P), and by
Proposition~\ref{p4.2}, the operator~$f_t$ is analytic with respect
to~$(u_0,\eta)$.  Furthermore, taking into account
Proposition~\ref{p4.5} and repeating the argument used in
Section~\ref{s3.1}, for any $u\in H$ and any positive $t\in J_T$ we
can find a ball~$B_u$ in a finite-dimensional subspace $X_u\subset
L^2(J_T,H)$ such that $f_t(u,B_u)\subset F$ has at least one interior
point. Thus, the conditions of Theorem~\ref{t2.2} are fulfilled, and
we can conclude that assertions~(i) and~(ii) hold.
\end{proof}

\begin{example}
Let us consider an example of a random force~$\eta(t,x)$ for which
the hypotheses of Theorem~\ref{t3.2} are satisfied for any~$T>0$. Let $H_0\subset H$
be any finite-dimensional subspace and let $\{\eta(t),t\ge0\}$ be a homogeneous Gaussian
process in~$H_0$ with a correlation function~$K(t)$. (For the existence of such
a process, see~\cite[Chapter~3]{GS}.) Suppose that~$K(t)$ is a positive-definite operator
for any $t\ge0$. Then, for any $T>0$, the distribution of the restriction of~$\eta$
to~$J_T$ is a non-degenerate Gaussian measure. As is explained in Example~\ref{e1.4},
such a measure satisfies property~(P).
\end{example}

\subsection{Navier--Stokes equations perturbed by a white noise}
\label{s3.3}
This subsection is devoted to studying the NS system~\eqref{3.1}, \eqref{3.15},
in which $h\in L_{\rm loc}^2(\R_+,H)$ is a deterministic function and~$\eta(t,x)$
is a random process white in time and $H^2$-regular in the space variables.
More precisely, we define the space $U=V\cap H^2(\T^2,\R^2)$ endowed with
the $H^2$-norm and assume that there is a non-negative nuclear operator
$Q\in\LL(U)$ such that
\begin{equation} \label{3.16}
\eta(t,x)=\frac{\p}{\p t}\zeta(t,x),
\end{equation}
where~$\zeta$ is a Gaussian process in~$U$ with continuous
trajectories and the covariance operator
$$
K(t,s)=(t\wedge s)Q, \quad t,s\ge0.
$$
It follows that
$$
\zeta(t,x)=\sum_{j=1}^\infty b_j\beta_j(t)g_j(x),
$$
where $\{g_j\}$ is an orthonormal basis in~$U$ formed of the
eigenvectors of~$Q$, $b_j^2$~is the eigenvalue of~$Q$ corresponding
to~$g_j$, and~$\{\beta_j\}$ is a sequence of independent standard
Brownian motions.  It is well known (see~\cite{VF,F94}) that for
almost every value of the random parameter the Cauchy problem
for~\eqref{3.1}, \eqref{3.15}, \eqref{3.16} has a unique solution in
the space~$\XX$, and we denote by $S_t:H\to H$ its resolving (random)
operator.

\begin{theorem} \label{t3.3}
There is an integer\,\footnote{A more precise description of~$N$ can
be found in Remark~\ref{r3.2}} $N\ge1$ not depending on~$\nu>0$ such
that if $h\in L_{\rm loc}^2(\R_+,H)$ is a given function and the image
of~$Q$ contains~$H_N$, then for any $u_0\in H$ and $t>0$
assertions~{\rm(i)} and~{\rm(ii)} of Theorem~\ref{t3.2} hold, and
\begin{equation} \label{3.17}
\supp\lambda_u(t)=F.
\end{equation}
\end{theorem}

\begin{proof}
{\it Step~1}.
Let us fix $T>0$ and study our problem on the interval~$J_T$.  The
solution of~\eqref{3.1}, \eqref{3.3}, \eqref{3.15}, \eqref{3.16} can
be represented in the form $u=v+\zeta$, where $v\in\XX_T$ is a
solution of the problem
\begin{equation} \label{3.18}
\dot v+\nu Lv+B(v,v)+B(v,\zeta)+B(\zeta,v)
=h(t)-\nu L\zeta-B(\zeta,\zeta), \quad v(0)=u_0.
\end{equation}
Let $\FF:H\times L^2(J_T,U)\to\XX_T$ be the operator that takes~$(u_0,\zeta)$
to the solution $v\in\XX_T$ of~\eqref{3.18}. Using standard methods of the theory
of 2D NS equations, one can show that~$\FF$ is well defined and analytic
(cf.~Proposition~\ref{p4.2}). Thus, denoting by~$\FF_t(u_0,\zeta)$ the restriction
of~$\FF(u_0,\zeta)$ to the time~$t$, we can write
$$
{\mathsf P}_Fu(t)={\mathsf P}_F\zeta(t)+{\mathsf P}_F\FF_t(u_0,\zeta)
=:f_t(u_0,\zeta).
$$
It follows that $\lambda_{u_0}(t)=\DD({\mathsf P}_Fu(t))$ is the image
of~$\mu_T$ under the analytic map~$f_t(u_0,\cdot)$, where~$\mu_T$ is
the distribution of the restriction of~$\zeta$ to~$J_T$. Hence,
assertions~(i) and~(ii) of Theorem~\ref{t3.2} will be established if
we show that the hypothesis of Theorem~\ref{t2.2} are satisfied.

\smallskip
{\it Step~2}.
Let us regard~$\mu_T$ as a measure on~$L^2(J_T,U_0)$, where~$U_0$
denote the closure in~$U$ of the image of~$Q$. In this case,
$\mu_T$~is a non-degenerate Gaussian measure that satisfies
property~(P) (see Examples~\ref{e1.2} and~\ref{e1.4}). We claim that
for any $u_0\in H$, $t\in J_T\setminus\{0\}$, and~$r>0$ there is a
ball~$B_{u_0}=B_{u_0}(t,r)$ in a finite-dimensional subspace
of~$L^2(J_T,H_N)$ such that
\begin{equation} \label{3.19}
f_t(u_0,B_{u_0})\supset B_F(r).
\end{equation}
Indeed, for any deterministic function $\xi\in C^1(J_T,H_N)$, we can
write
\begin{equation} \label{3.20}
u(t)=\xi(t)+\FF_t(u_0,\xi)=\RR_t(u_0,\hat\xi),
\end{equation}
where~$\RR_t$ is the resolving operator for the NS system (see
Section~\ref{s4.3}) and $\hat\xi=\p_t\xi$. It follows from
Proposition~\ref{p4.5} (cf.~Step~2 in the proof of Theorem~\ref{t3.1})
that if~$N\ge1$ is sufficiently large, then there is a ball~$\widehat
B_{u_0}$ in a finite-dimensional subspace of $C(J_T,H_N)$ such that
\begin{equation} \label{3.21}
{\mathsf P}_F\RR_t(u_0,\widehat B_{u_0})\supset B_F(r).
\end{equation}
It is clear that the image of~$\widehat B_{u_0}$ under the linear
operator $\hat\xi\mapsto\int_0^\cdot\hat\xi(s)\,ds$ is contained in a
finite-dimensional ball~$B_{u_0}$. Combining this with~\eqref{3.20}
and~\eqref{3.21}, we arrive at~\eqref{3.19}.  Hence,
Theorem~\ref{t2.2} is applicable, and we obtain assertions~(i)
and~(ii).

\smallskip
{\it Step~3}.
Relation~\eqref{3.17} is a standard consequence of controllability,
and we confine ourselves to outlining its proof. Let $y\in F$ be an
arbitrary point and~$\e>0$. We wish to show that
\begin{equation} \label{3.22}
\lambda_{u_0}(t)(B_F(y,\e))>0.
\end{equation}
In view of~\eqref{3.19}, there is $\xi_0\in L^2(J_T,H_N)$ such that
$f_t(u_0,\xi_0)=y$.  By continuity, there is $\delta>0$ such that if
\begin{equation} \label{3.23}
\|\xi-\xi_0\|_{L^2(J_T,U)}\le\delta,
\end{equation}
then $\|f_t(u_0,\xi)-f_t(u_0,\xi_0)\|_F\le\e$. Since~$\mu_T$ is non-degenerate,
it follows from~\eqref{1.7} that
$$
\mu_T\bigl(\{\xi\in L^2(J_T,U_0):\mbox{$\xi$ satisfies~\eqref{3.23}}\}\bigr)>0
\quad\mbox{for any $\delta>0$}.
$$
This implies the required property~\eqref{3.22}. The proof is complete.
\end{proof}

\begin{remark} \label{r3.4}
It is established in~\cite{AS-2006} that the 2D Euler equations
considered in the space~$H^s(\T^2,\R^2)$ with $s\ge2$ possesses the
property of exact controllability for observed projections
(cf.~Section~\ref{s4.3}). Using this fact, one can show that some
results similar to those in Theorems~\ref{t3.2} and~\ref{t3.3} are
true for the Euler equations. This will be done in a forthcoming
publication.
\end{remark}

As in the case of the NS system perturbed by the piecewise-constant
random force~\eqref{3.2}, Theorem~\ref{t3.3} readily implies the
existence of a density with respect to the Lebesgue measure for
finite-dimensional projections of any stationary distribution. Namely,
consider Eq.~\eqref{3.1} with the right-hand side
\begin{equation} \label{3.02}
g(t,x)=h(x)+\eta(t,x),
\end{equation}
where $h\in H$ is a deterministic function and~$\eta$ is a random
process of the form~\eqref{3.16}. In this case, problem~\eqref{3.1},
\eqref{3.02} generates a Markov family in the space~$H$, which has at
least one stationary measure (e.g., see~\cite{F94}). The proof of the
following result is similar to that of Corollary~\ref{c3.4} and will
be omitted.

\begin{corollary} \label{c3.9}
Under the conditions of Theorem~\ref{t3.3}, for any deterministic
function $h\in H$, any stationary measure~$\mu$, and any
finite-dimensional subspace~$F\subset H$, the 
projection~${\mathsf P}_F\mu$ possesses an almost everywhere positive
density with respect to the Lebesgue measure on~$F$.
\end{corollary}

\section{\sloppy Appendix}
\label{s4}

\subsection{Small perturbations of smooth functions}
\label{s4.1}
Let~$E$ and~$F$ be finite-dimensional vector spaces,
let $D\subset E$ be an open subset, and
let $f:D\to F$ be an infinitely differentiable function.

\begin{proposition} \label{p4.1}
Suppose that $f(D)$ has a non-empty interior. Then there is~$\e>0$ such
that if $g:D\to F$ is a continuous function satisfying the inequality
\begin{equation} \label{4.1}
\sup_{x\in D}\|f(x)-g(x)\|_F\le\e,
\end{equation}
then $g(D)$ has a non-empty interior.
\end{proposition}

\begin{proof}
Since the interior of~$f(D)$ is non-empty, Sard's theorem
(see~\cite{Sternberg}) implies that~$f$ has at least one regular
point~$x_0\in D$. Let us set $y_0=f(x_0)$ and assume,
without loss of generality, that $x_0=0$ and $y_0=0$.
We can represent~$E$ as a direct sum $E=E_1\dotplus E_2$ such that the
restriction of the derivative~$Df(x_0)$ to~$E_1$ is an isomorphism
from~$E_1$ onto $F$. The required assertion will be proved if we show
that for any continuous function $g:D\to F$ satisfying~\eqref{4.1} with~$\e\ll1$
the interior of~$g(D\cap E_1)$ is non-empty. To this end, denote
by~$f_1$ and~$g_1$ the restriction of~$f$ and~$g$ to~$D\cap E_1$ and note
that if~$\delta>0$ is sufficiently small, then $\deg_\delta f_1=1$,
where~$\deg_\delta h$ denotes the degree of a continuous function $h:D\cap E_1\to F$
with the respect to the ball~$B_{E_1}(\delta)$. It follows that if~$g$ satisfies~\eqref{4.1}
with sufficiently small~$\e$, then $\deg_\delta g_1=\deg_\delta f_1$. Since $y=0$
is a regular point of~$f_1$, we can find~$r>0$ such that
any point  $y\in F$ with norm not exceeding~$r$ has a
preimage under~$g_1$. This shows that the interior of~$g_1(D)$ is non-empty.
The proposition is proved.
\end{proof}

\subsection{Resolving operator for the Navier--Stokes equations}
\label{s4.2}
Let us consider the Navier--Stokes system~\eqref{3.1} supplemented
with the initial condition~\eqref{3.3}. It is well-known~\cite{CF}
that if $u_0\in H$ and $g\in L^1(J_T,H)$ for some~$T>0$, then
problem~\eqref{3.1}, \eqref{3.3} has a unique solution
$$
u_\nu\in \XX_T:=C(J_T,H)\cap L^2(J_T,V).
$$
Let $\R_+^*=(0,\infty)$ and let
$\RR:\R_+^*\times H\times L^1(J_T,H)\to\XX_T$ be an operator that
takes $(\nu,u_0,f)$ to the solution $u_\nu\in\XX_T$. The following
result on analyticity of~$\RR$ is established\footnotemark{} in~\cite{kuksin}.
\footnotetext{In~\cite{kuksin}, the analyticity of~$\RR$ is proved for
  a fixed~$\nu>0$. However, an unessential modification of the proof
enables one to handle the general case.}

\begin{proposition} \label{p4.2}
The operator~$\RR$ is analytic on the domain of its definition.
\end{proposition}

We now consider the problem~\eqref{3.1}~--~\eqref{3.3}, in which $T>0$
is a parameter and $\{\eta_k\}\subset H$ is a sequence of
deterministic functions. For any integer $k\ge1$, we set
$\HHH_k=H\times\cdots\times H\,\mbox{($k$ times)}$ and define the
operator
\begin{equation} \label{4.4}
\sS_k:\R_+^*\times H\times \HHH_k\to H, \quad
(T,u_0,\eeta_k)\to u(kT),
\end{equation}
where $\eeta_k=(\eta_1,\dots,\eta_k)\in\HHH_k$ and~$u(t)$ denotes the
solution of~\eqref{3.1}~--~\eqref{3.3}.

\begin{corollary} \label{c4.3}
For any integer~$k\ge1$, operator~\eqref{4.4} is analytic on the
domain of its definition.
\end{corollary}

\begin{proof}
Let us rewrite~\eqref{3.1} in the form
\begin{equation} \label{4.5}
\dot u+\nu Lu+B(u,u)=\sum_{k=1}^\infty I_{[k-1,k)}(t/T)\eta_k(x),
\end{equation}
where~$I_{[k-1,k)}$ stands for the indicator function
of~$[k-1,k)$. Performing the change of variables $u(t)=T^{-1}v(s)$,
$s=t/T$, in Eqs.~\eqref{4.5} and~\eqref{3.3}, we arrive at the problem
$$
\p_s v+T\nu Lv+B(v,v)=T^2\hat g(s,x),
\quad v(0)=Tu_0,
$$
where $\hat g(s,x)=\sum_k I_{[k-1,k)}(s)\eta_k(x)$. It follows that
$$
\sS_k(T,u_0,\eeta_k)=T^{-1}\RR_k(T\nu,Tu_0,T^2\hat g),
$$
where $\RR_t(\nu,u_0,f)$ denotes the restriction of~$\RR(\nu,u_0,f)$
to the time~$t$. The required result follows now from
Proposition~\ref{p4.2}.
\end{proof}

\subsection{A controllability property of the Navier--Stokes equations}
\label{s4.3}
Recall that $J_T=[0,T]$ and that~$H_N\subset H$ denotes the vector
space spanned by the first~$N$ eigenfunctions of the Stokes
operator~$L$. Consider the Navier--Stokes system~\eqref{3.1}, \eqref{3.15}
where $h\in L^2(J_T,H)$ is a given function and $\eta\in L^2(J_T,H_N)$
is a control function. In this subsection, we assume that the
viscosity $\nu>0$ and the function~$h$ are fixed and denote by
$\RR_t:H\times L^2(J_T,H)\to H$ the operator that takes~$(u_0,\eta)$
to the function~$u(t)$, where $u\in\XX_T$ is the solution
of~\eqref{3.1}, \eqref{3.3}, \eqref{3.15}.

Let $F\subset H$ be a finite-dimensional subspace,
let ${\mathsf P}_F:H\to F$ be the orthogonal projection in~$H$ onto~$F$,
and let~$T>0$ be a constant.

\begin{definition} \label{d4.4}
Problem~\eqref{3.1}, \eqref{3.15} is said to be {\it controllable in time~$T$
for the projection to~$F$\/} if ${\mathsf P}_F\RR_T(u_0,L^2(J_T,H_N))\supset F$
for any~$u_0\in H$.

Problem~\eqref{3.1}, \eqref{3.15} is said to be {\it solidly
controllable in time~$T$ for the projection to~$F$\/} if for any~$R>0$
and~$u_0\in H$ there is a constant~$\e>0$ and a compact subset
$K=K(R,u_0)\subset L^2(J_T,H_N)$ such that $\varPhi(K)\supset B_F(R)$,
where $\varPhi:K\to F$ is an arbitrary continuous map satisfying the
inequality
$$
\sup_{\eta\in K}\|\varPhi(\eta)-{\mathsf P}_F\RR_T(u_0,\eta)\|_F\le\e.
$$
\end{definition}

\begin{proposition} \label{p4.5}
There is an integer~$N\ge1$ not depending on~$\nu$ and~$h$ such that
for any~$T>0$ and any finite-dimensional subspace~$F\subset H$ the
Navier--Stokes system~\eqref{3.1}, \eqref{3.15} is solidly
controllable in time~$T$ for the projection to~$F$.
\end{proposition}

In the case $h\equiv0$, this result is established
in~\cite{AS-2005,AS-2006}. The general situation can be treated by the
same argument.

We now give a more precise description of the integer~$N$
entering Proposition~\ref{p4.5}. To this end, it is convenient
to index the trigonometric basis in~$H$ by the elements of~$\Z^2$.
Namely, we write $j=(j_1,j_2)\in\Z^2$ and set
\begin{gather*}
e_j(x)=\sin(jx)\,j^\bot\quad
\mbox{for $j_1>0$ or $j_1=0$, $j_2>0$},\\
e_j(x)=\cos(jx)\,j^\bot\quad
\mbox{for $j_1<0$ or $j_1=0$, $j_2<0$},\\
e_0^1(x)=(1,0), \quad e_0^2(x)=(0,1),
\end{gather*}
where $j^\bot=(-j_2,j_1)$.
The family
$\EE=\{e_0^i,e_j,\,i=1,2,j\in\Z^2\setminus\{0\}\}$
is a complete set of eigenfunctions for the Stokes
operator~$L$ and, hence, is an orthogonal basis in~$H$.

For any symmetric subset $\KK\subset\Z^2$ containing the
point~$(0,0)$, we write $\KK_0=\KK$ and define~$\KK^i$ with $i\ge1$ as
the union of~$\KK^{i-1}$ and the family of vectors $l\in\Z^2$ for
which there are $m,n\in\KK^{i-1}$ such that
$$
l=m+n,\quad |m|\ne|n|,\quad m\wedge n\ne0,
$$
where $m\wedge n=m_1n_2-m_2n_1$.

\begin{definition}
A symmetric subset~$\KK\subset\Z^2$ containing~$(0,0)$ is said to be
{\it saturating\/} if $\cup_{i\ge0}\KK^i=\Z^2$.
\end{definition}

For any subset $\KK\subset\Z^2$, we denote by $H_\KK$ the vector space
spanned by the family $\{e_0^1,e_0^2,e_j,j\in\KK\setminus0\}$. The
following result is a refined version of Proposition~\ref{p4.5}
(see~\cite{AS-2005,AS-2006}).

\begin{proposition}
Let $\KK\subset\Z^2$ be a saturating subset. Then for any positive
constant~$\nu$ and~$T$, any function $h\in L^2(J_T,H)$, and any
finite-dimensional subspace~$F\subset H$, the Navier--Stokes
system~\eqref{3.1}, \eqref{3.15} with $\eta\in L^2(J_T,H_\KK)$ is
solidly controllable in time~$T$ for the projection to~$F$.
\end{proposition}

\end{document}